\documentclass[11pt,a4paper,reqno]{amsart}

\title[$q$-Narayana numbers]{$q$-Narayana numbers and \\the flag $h$-vector of
$J({\bf 2} \times {\bf n})$}

\author{Petter Bränd\'{e}n}
\address{Matematik,
  Chalmers tekniska h\"ogskola och G\"oteborgs universitet,\linebreak
  S-412~96  G\"oteborg, Sweden}
\email{branden@math.chalmers.se}
\keywords{Narayana numbers, flag h-vector, Schur Function, shelling}
\date{\today}

\usepackage[english]{babel}
\usepackage{amsmath}
\usepackage{amsthm,amssymb,latexsym}
\usepackage{t1enc}
\usepackage{epic}
\usepackage{eepic}
\usepackage{xypic}
\xyoption{all}

\newtheorem{prop}{Proposition}
\newtheorem{lemma}[prop]{Lemma}
\newtheorem{cor}[prop]{Corollary}
\newtheorem{thm}[prop]{Theorem}

\theoremstyle{definition}

\newtheorem{example}[prop]{Example}

\newcommand{\N}{\mathbb{N}}
\newcommand{\Z}{\mathbb{Z}}
\newcommand{\JH}{\mathcal{L}}
\newcommand{\la}{<^{{\tiny \Lambda}}}

\newcommand{\lom}{<^{{\tiny \Omega}}}

\newcommand{\lomq}{\leq^{{\tiny \Omega}}}
\newcommand{\Long}{ LS }
\newcommand{\Dy}{\mathcal{D}}
\newcommand{\one}{\hat{1}}
\newcommand{\zero}{\hat{0}}
\newcommand{\Max}{\mathfrak{M}}
\newcommand{\e}{\ar@{-}}
\newcommand{\de}{\ar@{.}}
\def\gen#1,{\langle #1 \rangle}
\def\geno#1,{\langle #1 \rangle_{\infty}}
\def\qbin#1,#2,{{#1\atopwithdelims[]#2}}
\def\newop#1{\expandafter\def\csname #1\endcsname{\mathop{\rm
#1}\nolimits}}
\newop{row}
\newop{des}
\newop{hp}
\newop{da}
\newop{ea}
\newop{lnfs}
\newop{MAJ}
\newcommand{\Discup}[2]{\bigcup_{#1}^{#2}\hspace{-2ex}\cdot\hspace{1ex}}
\newcommand{\discup}{\makebox{$\cup\hspace{-1.1ex}\cdot\hspace{.55ex}$}}

\usepackage{graphicx}

\begin{document}
\maketitle
\thispagestyle{empty}
\begin{abstract}
  The Narayana numbers are $N(n,k) = {1 \over n}{n \choose k}{n
    \choose {k+1}}$.  There are several natural statistics on Dyck
  paths with a distribution given by $N(n,k)$.  We show the
  equidistribution of Narayana statistics by computing the flag
  $h$-vector of $J({\bf 2} \times {\bf n})$ in different ways. In the
  process we discover new Narayana statistics and provide
  co-statistics for which the Narayana statistics in question have a
  distribution given by F\"urlinger and Hofbauers $q$-Narayana
  numbers. We  also interpret the $h$-vector in terms of 
  semi-standard Young tableaux, which enables us to express the
  $q$-Narayana numbers in terms of Schur functions.
\end{abstract}  

\section{Introduction}
The {\em Narayana numbers}, 
$$
N(n,k) = {1 \over n}{n \choose k}{n \choose {k+1}},
$$
appear in many combinatorial problems. Some examples are the 
number of noncrossing partitions of $\{1,2, \ldots, n\}$ of rank 
$k$ \cite{krew3}, the number of $132$-avoiding permutations with 
$k$ descents \cite{simion}, and also several problems involving 
Dyck paths. 

A {\em Dyck path} of length $2n$ is a path in $\N \times \N$ 
from $(0,0)$ to $(n,n)$ using steps $v=(0,1)$ and $h=(1,0)$, which never 
goes below the line $x=y$. The set of all Dyck paths of length $2n$ is 
denoted $\Dy_n$. A statistic on $\Dy_n$ having a 
distribution given by the Narayana numbers
will in the sequel be  referred to as a {\em Narayana statistic}.
The first Narayana statistics to be discovered were
\begin{itemize}
\item[$\des(w)$:] the number of {\em descents (valleys)} (sequences
  $hv$) in $w$, \cite{narayana1},
 \item[$\ea(w)$:] the number of {\em even ascents}, i.e., the number of 
letters $v$ in an even position in $w$, \cite{krew2}, 
 \item[$\lnfs(w)$:] the number of {\em long non-final sequences}, more 
       precisely the number of sequences $vvh$ and $hhv$ in $w$, 
       \cite{krew1}. 
\end{itemize} 
Recently, \cite{deutsch}, a new Narayana statistic, $\hp$, was
discovered and it counts the number of {\em high peaks}, i.e., peaks
not on the diagonal $x=y$.  Also, in \cite{sulanke1,sulanke2}
Sulanke found numerous new Narayana statistics with the help of a
computer. For terminology on posets in what follows, we refer the reader 
to \cite{stanley1}.

We will show that $\des, \hp$ and $\lnfs$ arise when computing the flag
$h$-vector of the  lattice $J( {\bf 2} \times {\bf n})$ of order
ideals in the poset ${\bf 2} \times {\bf n}$ in different ways.
In Section \ref{deshp} we will show how the statistics descents and 
high peaks arise when considering different linear extensions of 
${\bf 2} \times {\bf n}$. This will give the equidistribution of the 
descent-set and the set of high-peaks. In Section \ref{lnfss} we consider 
a shelling of the order complex $\Delta(J( {\bf 2} \times {\bf n}))$ to 
show that the set of long non-final sequences has the same distribution as 
the descent set over Dyck paths.

There is a $q$-analog of the Narayana numbers,
$$
N_q(n,k)=  \frac {1} {[n]} \qbin n,k, \qbin n, {k+1},q^{k^2+k},
$$
introduced by F\"{u}rlinger and Hofbauer in \cite{furhof}.  To each
statistic we treat we will associate a co-statistic  together with
which the Narayana statistic has a joint distribution given by the
$q$-Narayana numbers.

\section{Descents and High peaks}\label{deshp}
Let $P$ be any finite graded poset with a smallest element $\zero$ and a 
greatest element $\one$ and let $\rho$ be the rank function of $P$ with   
$\rho(P) := \rho(\one)=n$. For $S \subseteq [n-1]$ let    
$$
\alpha_P(S) := |\{ c \mbox{ is a chain of } P : \rho(c)=S\}|,
$$
and 
$$
\beta_P(S):=\sum_{T \subseteq S} (-1)^{|S-T|}\alpha_P(T).  
$$
The functions $\alpha_P, \beta_P :  2^{[n-1]} \rightarrow \Z$ are 
the {\em flag f-vector} and the {\em flag h-vector} of $P$ respectively.

If  $P$ is a finite poset of cardinality $p$ and  
$\omega : P \rightarrow [p]$ is a linear extension of $P$ then the 
{\em Jordan-Hölder set}, $\JH(P, \omega)$, of $(P, \omega)$
is the set of permutations $a_1a_2 \cdots a_p$ such that 
$\omega^{-1}(a_1),\omega^{-1}(a_2), \ldots, \omega^{-1}(a_p)$ is a 
linear extension of $P$, in other words 
$$
\JH(P, \omega) = \{ \omega \circ \sigma^{-1} : \sigma 
\mbox{ is a linear extension of } P \}. 
$$
We will need the following theorem (Theorem 3.12.1 of \cite{stanley1}):
\begin{thm}\label{stanley}
Let $L=J(P)$ be a distributive lattice of rank $p=|P|$, and let $\omega$ be a 
linear extension of $P$. Then for all $S \subseteq [p-1]$ we have that  
$\beta_L(S)$ is equal to the number of permutations $\pi \in \JH(P, \omega)$ 
with descent set $S$.
\end{thm} 

It will be convenient to code a Dyck path $w$ in the letters
$\{v_i\}_{i=1}^\infty \cup \{h_i\}_{i=1}^\infty$ by letting $v_i$ and
$h_i$ stand for the $i$th vertical step and the $i$th horizontal step
 in $w$, respectively. Thus $vvhvhh$ is coded as
$v_1v_2h_1v_3h_2h_3$.  We may write  the set of elements of ${\bf
  2} \times {\bf n}$ as the disjoint union $C_1 \discup C_2$ where
$C_i = \{ (i,k): k \in [n] \}$ for $i=1,2$. For any linear extension
$\sigma$ of ${\bf 2} \times {\bf n}$ let $W(\sigma)$ be the Dyck path
$w_1w_2 \cdots w_{2n}$ where
$$ 
w_i=
\begin{cases}
v_j \ \ \mbox{ if } \sigma^{-1}(i) = (1,j) \mbox{ and } \\
h_j \ \ \mbox{ if } \sigma^{-1}(i) = (2,j). 
\end{cases}
$$
It is 
clear that $W$ is a bijection between the set of linear extensions 
of ${\bf 2}\times {\bf n}$ and the set of Dyck paths of length $2n$. 

\newcommand{\n}{*=0{\bullet}}
\begin{center}
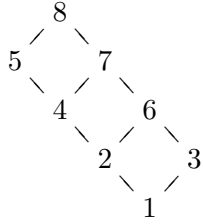
\begin{figure}\caption{\label{fig3} The linear extension of 
${\bf 2} \times {\bf 4}$ corresponding to the Dyck path 
$v_1v_2h_1v_3v_4h_2h_3h_4$.}
$$
\vcenter{\xymatrix@R=5pt@C=5pt{
          & 8           &                &                &   \\
 5\e[ru]&                & 7\e[lu]     &                &   \\
          &4\e[lu]\e[ru]&                &  6\e[lu]    &    \\
          &                &2\e[lu]\e[ru]&                & 3\e[lu] \\
          &                &                &1\e[lu]\e[ru]& 
    }}
$$
\end{figure}
\end{center}
Fix a Dyck path $W_0 \in \Dy_n$ and let $\omega_0=W^{-1}(W_0)$. Now, 
if 
$\pi = \omega_0 \circ \sigma^{-1} \in \JH({\bf 2} \times {\bf n},\omega_0)$  
let $W(\sigma)= w_1w_2 \cdots w_{2n}$. 
Then $\pi(i) > \pi(i+1)$ if and only if $w_{i+1}$ comes before $w_i$ in 
$W_0$.   
%
%
In light of this we define, given Dyck paths $W_0$ and 
$w=w_1w_2 \cdots w_{2n}$,   
the {\em descent set of $w$ with respect to $W_0$} as 
$$
D_{W_0}(w) = \{ i \in [2n-1]:  w_{i+1} \mbox{ comes before } w_i \mbox{ in }
                W_0 \}.
$$
The descent set of $v_1h_1v_2v_3h_2h_3$ with respect to $v_1v_2h_1v_3h_2h_3$ 
is thus $\{ 2 \}$.
By  Theorem \ref{stanley} we now have:
\begin{thm}\label{main}
Let $W$ be any Dyck path of length $2n$ and let $S \subseteq
[2n-1]$ and let $\beta_n = \beta_{J({\bf 2} \times {\bf n})}$.
Then  
$$
\beta_n(S)=|\{ w \in \Dy_n : D_{W}(w) = S \}|.
$$
\end{thm}
For a 
given Dyck path $W$ we define the statistics $\des_{W}$, and 
$\MAJ_{W}$ by 
\begin{eqnarray*}
\des_{W}(w) &=& |D_{W}(w)|,\\
\MAJ_{W}(w) &=& \sum_{i \in D_{W}(w)}i .
\end{eqnarray*}
\begin{example}Two known Narayana statistics arise when fixing $W$ in certain 
ways: 
\begin{itemize}
\item[a)]
If $W=v_1v_2 \cdots v_n h_1h_2 \cdots h_n$ then $\des_{W}=\des$. 
\item[b)]
If $W=v_1h_1v_2h_2 \cdots v_nh_n$ then $\des_{W}=\hp$. Thus as a 
consequence of Theorem \ref{main} we have that the number of valleys and the 
number of high peaks have the same distribution over $\Dy_n$. This was 
first proved by Deutsch in \cite{deutsch}.  
\item[c)]
If $W = v_1h_1v_2v_3\cdots v_nh_2h_3\cdots h_n$ then $\des_{W}$ counts 
valleys $h_iv_j$ where $i>1$ and high peaks of the form $v_ih_1$.
\end{itemize}
\end{example} 
When $W=v_1v_2 \cdots v_n h_1h_2 \cdots h_n$ we 
drop the subscript and let $\des = \des_{W}$ and $\MAJ=\MAJ_{W}$. 
In \cite{furhof} F\"{u}rlinger and Hofbauer defined the 
$q$-{\em Narayana numbers}, $N_q(n,k)$, by  
$$
N_q(n,k):=\sum_{w \in \Dy_n, \des(w)=k} q^{\MAJ(w)}.
$$
We say that the bi-statistic $(\des, \MAJ)$ has the $q$-Narayana 
distribution. We will later see that $N_q(n,k)$ can be written in an 
explicit form.    
By Theorem~ \ref{main} we now have:
\begin{cor}
For all $W \in \Dy_n$ the bi-statistic $(\des_{W},\MAJ_{W})$ has 
the $q$-Narayana distribution.
\end{cor}
Let $\lambda=(\lambda_1, \lambda_2, \ldots, \lambda_{\ell})$ be a partition 
of positive integers. The index $\ell$ is called the {\em length}, 
$\ell(\lambda)$, of $\lambda$. A 
{\em semistandard Young tableau} (SSYT) of {\em shape} $\lambda$ is an 
array $T = (T_{ij})$ of positive integers, where 
$1 \leq i \leq \ell(\lambda)$ and 
$1 \leq j \leq \lambda_i$, that is weakly increasing in every row and 
strictly increasing in every column.   
\begin{figure}\caption{ An example of a SSYT of shape 
$(6,5,4,4,2)$ .}
$$
\begin{array}{ccccccccccc}
1&2&2&3&5&5\\
2&3&4&4&6&&\\      
4&5&5&6&&&&\\ 
5&6&6&8&&&&\\
7&8&&&&.&&&
\end{array}
$$
\end{figure}
For any SSYT of shape $\lambda$ let 
$$
x^T := x_1^{\alpha_1(T)}x_2^{\alpha_2(T)}\cdots,
$$
where $\alpha_i(T)$ denotes the number of entries of $T$ that are equal to 
$i$. The {\em Schur function} $s_{\lambda}(x)$ of {\em shape} $\lambda$ is the 
formal power series
$$
s_{\lambda}(x)= \sum_T x^T,
$$
where the sum is over all SSYTs $T$ of shape $\lambda$.  
If $T$ is any SSYT we let 
$\row(T)=(\gamma_1(T), \gamma_2(T), \ldots )$ where 
$\gamma_i(T) = \sum_j T_{ij}$. Let $\langle 2^k\rangle$ be the 
partition $(2,2, \ldots, 2)$ with $k$ 2's. 
\begin{thm}\label{ssyt}
For any $n > 0$ and $S \subseteq [2n-1]$, $|S|=k$, we have that
$\beta_n(S)$ counts 
the number of {\em SSYT}s $T$ of shape $\langle 2^k\rangle $ with $\row(T)=S$ and with 
parts less than $n$.  
\end{thm} 
\begin{proof}
Let $T$ be a SSYT as in the statement of the theorem. We want to 
construct a Dyck path $w(T)$ with descent set $S$. 

Let 
$w(T)=w_1w'_1w_2w'_2 \cdots w_{k+1}w'_{k+1}$ where 
\begin{itemize}
\item $w_1$ is the word  consisting of $T_{12}$ vertical steps and
  $w_1'$ is the word   consisting of $T_{11}$ horizontal steps,
\item $w_i$ is the word  consisting of $T_{i2}-T_{(i-1)2}$ vertical
  steps and $w_i'$ is the word  consisting of $T_{i1}-T_{(i-1)1}$
  horizontal steps, when $2 \leq i \leq k$,
\item $w_{k+1}$ is the word   consisting of $n-T_{k2}$ vertical steps and 
      $w_{k+1}'$ is the word   consisting of $n-T_{k1}$ horizontal steps.
\end{itemize}
It is clear that $w(T)$ is indeed a Dyck path with descent set $S$, and 
 each  such Dyck path is given by $w(T)$ for a unique SSYT $T$.
\end{proof}
\begin{center}
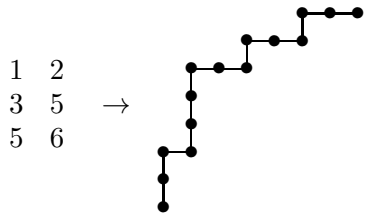
\begin{figure}\caption{\label{fig0} An illustration of Theorem \ref{ssyt} 
for $n=7$.}
$$
\begin{array}{cc}
1 & 2 \\
3 & 5 \\
5 & 6 
\end{array}
\;\;\rightarrow\;\;
\vcenter{\xymatrix@R=4pt@C=4pt@!{
&&&&&\n\e[r]&\n\e[r]&\n\\
&&&\n\e[r]&\n\e[r]&\n\e[u]\\
&\n\e[r]&\n\e[r]&\n\e[u]\\
&\n\e[u]\\
&\n\e[u]\\
\n\e[r]&\n\e[u]\\  
\n\e[u]\\
\n\e[u]
    }}
$$
\end{figure}
\end{center}
\begin{thm}\label{schur}
For all $n,k \geq 0$ we have 
$$
N_q(n,k)= s_{\langle 2^k\rangle }(q,q^2, \ldots, q^{n-1}).
$$
\end{thm}
\begin{proof} 
By  Theorem \ref{ssyt} we have that 
\begin{eqnarray*}
\sum_{w \in \Dy_n,\des(w)=k}q^{\MAJ(w)} &=& \sum_{|S|=k}\beta_n(S)
q^{\sum_{s \in S}s}\\
&=&\sum_T q^{\sum T_{ij}},
\end{eqnarray*}
where the last sum is over all $SSYT$s $T$ of shape $\langle 2^k\rangle $ with 
parts less than $n$. By the combinatorial definition of the Schur function 
this is equal to 
$s_{\langle 2^k\rangle }(q,q^2, \ldots, q^{n-1})$, and the theorem follows.
\end{proof}
%
%
If we identify a partition $\lambda$ 
with its diagram $\{(i,j): 1 \leq j \leq \lambda_i\}$ then  
the {\em hook length}, $h(u)$, at $u=(x,y) \in \lambda$ is defined by 
$$
h(u) = |\{ (x,j) \in \lambda : j \geq y\}|+ 
       |\{ (i,y) \in \lambda : i \geq x\}|-1,
$$
and the {\em content}, $c(u)$, is defined by 
$$
c(u) = y-x.
$$ 
We will use a result on Schur polynomials, commonly 
referred to as the {\em hook-content formula},  
see \cite[Theorem 7.21.2]{stanley2}. Let 
$[n] := 1 + q + \cdots +q^{n-1}$, $[n]! := [n][n-1]\cdots [1]$ and 
$$
\qbin n, k, := \frac {[n]!}{[n-k]![k]!}.
$$
\begin{thm}[Hook-content formula]\label{schurlemma}
For any partition $\lambda$ and $n >0$, 
$$
s_{\lambda}(q,q^2,\ldots,q^n) = 
q^{\sum i\lambda_i}\prod_{u \in \lambda}\frac {[n+c(u)]}{[h(u)]}.
$$
\end{thm}
We now have an alternative proof of the following result which was  
proved in \cite{furhof}, and is a special case of a result of MacMahon, 
stated without proof in \cite[p. 1429]{macmahon}.
\begin{cor}[F\"{u}rlinger, Hofbauer, MacMahon]
The $q$-Narayana numbers are given by:
$$
N_q(n,k)=  \frac {1} {[n]} \qbin n,k, \qbin n, {k+1},q^{k^2+k} 
$$
\end{cor}
\begin{proof}
The Corollary follows from Theorem \ref{schur} after an elementary 
application of the hook-content formula, 
which is left to the reader.
\end{proof}
\section{Long Non-final Sequences}\label{lnfss}
In \cite{krew1} Kreweras and Moszkowski defined a new Narayana statistic, 
$\lnfs$. Recall that a  
{\em long non-final sequence} in a Dyck path is a subsequence 
of type $vvh$ or $hhv$, and that the statistic $\lnfs$ is defined as the 
number of long non-final sequences in the Dyck path. We define the 
{\em long non-final sequence set}, $\Long(w)$, of a  Dyck path 
$w =a_1a_2 \cdots a_{2n}$ to be 
$$
\Long(w) = \{ i \in [2n-1] : a_{i-1}a_ia_{i+1} = vvh \mbox{ or } 
                             a_{i-1}a_ia_{i+1} = hhv \}.
$$
We will show that 
$$
\beta_n(S) = |\{ w \in \Dy_n : \Long(w) = S \}|. 
$$
To prove this we need some definitions.

An (abstract) 
{\em simplicial complex} $\Delta$ on a vertex set $V$ is a 
collection of subsets $F$ of $V$ satisfying:
\begin{enumerate}
  \item[(i)] if $x\in V$ then $\{x\}\in\Delta$,
  \item[(ii)] if $F\in\Delta$ and $G\subseteq F$, then $G\in\Delta$.
\end{enumerate}
The elements of $\Delta$ are called {\em faces} and 
a maximal face (with respect to inclusion) is called a {\em facet}. 
A simplicial complex is said
to be {\em pure} if all its facets have the same cardinality. 
A linear partial order $\Omega$ on the set of facets of a pure 
simplicial complex $\Delta$ is a {\em shelling} if 
%
%
whenever $F \lom G$ there is an $x \in G$ and $E \lom G$ such that 
$$
F \cap G \subseteq E \cap G = G \setminus \{x\}.
$$
A simplicial complex which allows a shelling is said to be {\em shellable}. 
Instead of finding a particular shelling we will find a partial 
order on the set of facets with the property that every linear extension is a 
shelling. In our attempts to prove that our partial order had this property 
we found ourselves proving 
Theorem \ref{preshelling} and Corollary \ref{postpreshelling}. We 
therefore take the opportunity to take a general approach and 
define what we call a {\em pre-shelling}. Though we have found examples 
of pre-shellings implicit in the literature we have not found explicit 
references, so we will provide proofs.

Let $\Omega$ be a partial 
order on the set of facets of a pure simplical complex $\Delta$. The 
{\em restriction}, $r_{\Omega}(F)$, of a facet $F$ is the set 
$$
r_{\Omega}(F) = \{ x \in F : \exists E \mbox{ s.t }  E \lom F \mbox{ and } 
E \cap F = F \setminus \{x\} \ \ \}.
$$
We say that $\Omega$ is a {\em pre-shelling} if any of the equivalent
conditions in Theorem \ref{preshelling} are satisfied.
\begin{thm}\label{preshelling}
Let $\Omega$ be a partial order on the set of facets of a pure simplicial 
complex $\Delta$. Then the following conditions on $\Omega$ are equivalent:
\begin{enumerate}
\item[(i)] For all facets $F,G$ we have
$$
r_{\Omega}(F) \subseteq G \mbox{ and } r_{\Omega}(G) \subseteq F \ \ 
\Longrightarrow \ \ F=G.
$$
\item[(ii)] $\Delta$ is the disjoint union 
$$
\Delta = \Discup{F}{} [r_{\Omega}(F), F].
$$
\item[(iii)] For all facets $F,G$ 
$$
r_{\Omega}(F) \subseteq G \ \ \Rightarrow \ \ F \lomq G.
$$

\item[(iv)] For all facets $F,G$: if $F \ngeq^{ {\tiny \Omega}} G$ then 
there is an $x \in G$ and $E \lom G$ such that
$$
F \cap G \subseteq E \cap G = G \setminus \{x\}.
$$
\end{enumerate}
\end{thm}
\begin{proof}
(i) $\Rightarrow$ (ii): Let $F$ and $G$ be facets of $\Delta$. If there 
is an $H \in [r_{\Omega}(F), F]\cap [r_{\Omega}(G), G]$ then 
$r_{\Omega}(F) \subseteq G$ and $r_{\Omega}(G) \subseteq F$, so by 
(i) we have $F=G$. Hence the union is disjoint. Suppose that $H \in \Delta$, 
and let $F_0$ be a minimal element, with respect to $\Omega$, of the set 
$$
 \{ F: F \mbox{~is a facet and~}    H \subseteq F \}.
$$
If $r_{\Omega}(F_0) \nsubseteq H$ then let 
$x \in r_{\Omega}(F_0) \setminus H$ and  let $E \lom F_0$ be such that 
$F_0 \cap E = F_0 \setminus \{x\}$. Then $H \subseteq E$, contradicting 
the minimality of $F_0$. This means that $H \in [r(F_0), F_0]$.    

(ii) $\Rightarrow$ (i): If $r_{\Omega}(F) \subseteq G$ and 
$r_{\Omega}(G) \subseteq F$ we have that 
$F \cap G \in [r_{\Omega}(F),F]\cap [r_{\Omega}(G),G]$, which by 
(ii) gives us $F=G$.

(i) $\Rightarrow$ (iii): If $r_{\Omega}(F) \subseteq G$ then by (i) we
have either $F=G$ or $r_{\Omega}(G) \nsubseteq F$. If $F= G$ we have
nothing to prove, so we may assume that there is an $x \in
r_{\Omega}(G) \setminus F$.  Then, by assumption, there is a facet
$E_1 \lom G$ such that
$$
r_{\Omega}(F) \subseteq G \cap E_1 = G \setminus \{x\} \subset E_1. 
$$
If $E_1 =F$ we are done. Otherwise we continue until we get 
$$
F=E_k \lom E_{k-1} \lom \cdots \lom E_1 \lom G,
$$
and we are done. 

(iii) $\Leftrightarrow$ (iv): It is easy to see that (iv) is just the 
contrapositive of (iii)

(iii) $\Rightarrow$ (i): Immediate.
\end{proof}
The set of all partial orders on the same set is 
partially ordered by inclusion, i.e $\Omega \subseteq \Lambda$ if 
$x \lom y$ implies $x \la y$.   
\begin{cor}\label{postpreshelling}
Let $\Delta$ be a pure simplicial complex. Then 
\begin{enumerate}
\item[(i)] all shellings of $\Delta$ are pre-shellings, 
\item[(ii)] if $\Omega$ is a pre-shelling of $\Delta$ and $\Lambda$ is a 
partial order such that $\Omega \subseteq \Lambda$, then $\Lambda$ is 
a pre-shelling of $\Delta$ with $r_{\Lambda}(F) = r_{\Omega}(F)$ for all 
facets $F$. In particular,  the set of all pre-shellings of $\Delta$ is 
an upper ideal of the poset of all partial orders on the set of facets 
of $\Delta$,
\item[(iii)] all linear extensions of a pre-shelling are shellings, with 
the same restriction function.  
\end{enumerate} 
\end{cor}
\begin{proof}
(i): Follows immediately from Theorem \ref{preshelling}(iv). 

(ii): That $\Lambda$ is a pre-shelling follows from 
Theorem \ref{preshelling}(iv). If $F$ is a facet 
then by definition 
$r_{\Omega}(F) \subseteq r_{\Lambda}(F)$, and if 
$r_{\Omega}(F) \subset r_{\Lambda}(F)$ for some facet $F$ we would have  
a contradiction by Theorem \ref{preshelling}(ii).

(iii): Is implied by (ii).
\end{proof}

Let $P$ be a finite graded poset with a smallest element $\zero$ and a 
greatest element $\one$. The {\em order complex}, $\Delta(P)$, of $P$ 
is the simplicial complex of all chains of $P$. A simplicial complex 
$\Delta$ is {\em partitionable} if 
it can be written as 
\begin{equation}\label{partitioning}
\Delta = [r(F_1), F_1] \discup [r(F_2), F_2] \discup \cdots \discup 
[r(F_n),F_n],
\end{equation}
where each $F_i$ is a facet of $\Delta$ and $r$ is any function on the set 
of facets such that $r(F) \subseteq F$ for all facets $F$. The right hand 
side of 
($\ref{partitioning}$) is a {\em partitioning} of~ $\Delta$. By 
Theorem \ref{preshelling}(iii) we see that shellable complexes 
are partitionable. We need the following well known fact about partitionable 
order complexes. Let $\Max(P)$ be the set of maximal chains of $P$. 
\begin{lemma}\label{parth}
Let $\Delta(P)$ be partitionable and let 
\begin{equation}\label{prt2}
\Delta(P) = \Discup{c}{}[r(c), c]
\end{equation} 
be a partitioning of $\Delta(P)$. Then the flag $h$-vector is given by
$$
\beta_P(S) = |\{ c \in \Max(P) : \rho(r(c))=S \}|.
$$
\end{lemma}
\begin{proof}
Let 
$\gamma_P(S) = |\{ c \in \Max(P) : \rho(r(c))=S \}|$. Note that if $c$ is 
a maximal chain then $\rho(c) = [0,\rho(\one)]$.  
By (\ref{prt2}) we have 
\begin{eqnarray*}
\alpha_P(S) &=& |\{ c \in \Delta(P) : \rho(c) = S\}| \\
            &=& |\{ c \in \Max(P) : \rho(r(c)) \subseteq S \}| \\
            &=& \sum_{T \subseteq S} \gamma_P(T),
\end{eqnarray*}
which, by inclusion-exclusion, gives $\gamma_P(S)= \beta_P(S)$.
\end{proof} 
We will identify the set of facets of $\Delta(J(\bf{2} \times
\bf{n}))$ with $\Dy_n$, the set of Dyck paths of length $2n$. We
therefore seek a partial order on $\Dy_n$ which is a pre-shelling.
Let $S=S(\Dy_n)$ be the set of mappings with elements
$$ s_i(w) = \begin{cases} a_1 \cdots a_{i-1}vhva_{i+3}\cdots a_{2n} & 
                 \mbox{ if } a_ia_{i+1}a_{i+2} = vvh,\\
                a_1 \cdots a_{i-1}hvha_{i+3}\cdots a_{2n}  
                & \mbox{ if } a_ia_{i+1}a_{i+2} = hhv,\\
                w & \mbox{ otherwise,}
                \end{cases}$$
for $1 \leq i \leq 2n-2$. Define a relation $\Omega_n$, by $u \lom w$ whenever  
$u \neq w$ and $u = \sigma_1\sigma_2 \cdots \sigma_k(w)$ for some mappings 
$\sigma_i \in S$ (see Figure \ref{fig1}).
\begin{lemma}
The relation $\Omega_n$ on $\Dy_n$ is a partial order.
\end{lemma}
\begin{proof}
We need to prove that $\Omega_n$ is anti-symmetric. To do this we define   
a mapping $\sigma :\Dy_n \rightarrow \N \times \N$, where 
$\N \times \N$ is ordered lexicographically, with the property 
$$
u \lom w \Rightarrow \sigma(u) < \sigma(w).
$$ 
Define $\sigma(w) = (\da(w), \MAJ(w))$, where $\da(w)$ is the number of 
double ascents (sequences $vv$) in $w$. Now, suppose that $s_i \in S$ and 
$s_i(w) \neq w=a_1a_2 \cdots a_{2n}$.  
Then $\da(s_i(w)) \leq \da(w)$, and if we have equality we must 
have $a_{i-1}a_ia_{i+1}a_{i+2}= vvhv$ or 
     $a_{i-1}a_ia_{i+1}a_{i+2}= hhvh$ which implies 
$\MAJ(s_i(w))<\MAJ(w)$, so $\sigma$ has the desired properties.  
\end{proof}

\begin{center}
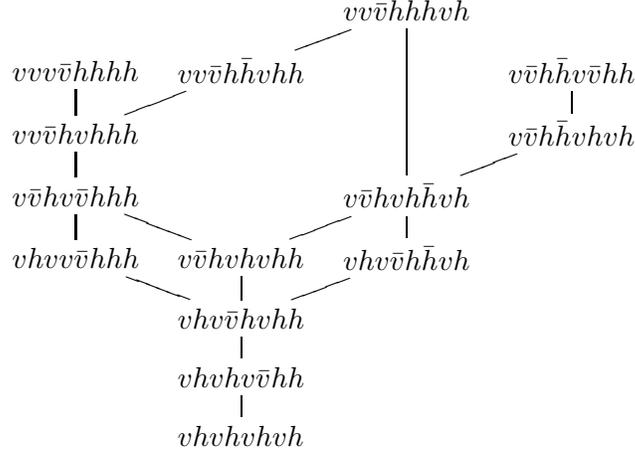
\begin{figure}
\caption{\label{fig1} The partial order $\Omega_4$ on $\Dy_4$, with long 
non-final sequences marked with bars.}
$$
\vcenter{\xymatrix@R=8pt@C=8pt{
             &                     & vv\bar{v}hh\bar{h}vh              &  \\
vvv\bar{v}hhhh  & vv\bar{v}h\bar{h}vhh\e[ru] &  
& v\bar{v}h\bar{h}v\bar{v}hh      \\
vv\bar{v}hvhhh\e[u]\e[ru] &  &    & v\bar{v}h\bar{h}vhvh\e[u]  \\ 
v\bar{v}hv\bar{v}hhh\e[u] &        & v\bar{v}hvh\bar{h}vh\e[uuu]\e[ru] &     \\ 
vhvv\bar{v}hhh\e[u] & v\bar{v}hvhvhh\e[lu]\e[ru] & vhv\bar{v}h\bar{h}vh\e[u] & \\ 
             &  vhv\bar{v}hvhh\e[u]\e[lu]\e[ru] &                 &     \\  
             &  vhvhv\bar{v}hh\e[u]             &                 &      \\  
             &  vhvhvhvh\e[u]             &                 &
    }}
$$
\end{figure}
\end{center}

If $v$ and $w$ intersect maximally then it is plain to see that 
either $v = s(w)$ or $s(v) = w$ for some $s \in S$. It 
follows that if $w = a_1a_2 \cdots a_{2n}$ then 
$$
r_{\Omega_n}(w) = \{ a_1+a_2+ \cdots + a_i : i \in \Long(w)\},
$$
so $\rho(r_{\Omega_n}(w))=\Long(w)$. It remains to prove that $\Omega_n$ is a 
pre-shelling.
\begin{thm}\label{longshelling}
For all $n \geq 1$ the partial order $\Omega_n$ is 
a pre-shelling of $\Dy_n$. 
\end{thm}
\begin{proof}
We prove that $\Omega_n$ satisfies the contrapositive of condition (i) of 
 Theorem \ref{preshelling}.
Suppose that $u=a_1a_2\cdots a_{2n} \neq w=b_1b_2\cdots b_{2n}$ and let 
$k$ be the coordinate such that $a_i = b_i$ for $i< k$ and 
$a_k \neq b_k$. By symmetry we may assume that $a_k=h$. Now, if $a_{k-1}=h$ 
then the valley of $u$ which is determined by the first $v$ 
(at, say, coordinate $\ell +1$) after  
$k$ will correspond to an element  
$$
x = a_1 + \cdots + a_{\ell} \in r_{\Omega_n}(u)\setminus w
$$
(see Figure \ref{fig2}). 
\begin{center}
\begin{figure}\caption{\label{fig2}}
$$
\vcenter{\xymatrix@R=12pt@C=12pt{
          &                     &         &             &    \\
          &    \n\de[u]         &         &             & \n   \\
 \n\e[r]_{a_{k-1}}  & \n\e[r]_{a_{k}}\de_{b_k}[u]       & \n\e[r] & 
\ldots\e[r]_{\ \ a_{\ell}} & \n \e_{a_{\ell+1}}[u]   
    }}
$$
\end{figure}
\end{center}
If $a_{k-1}=v=b_{k-1}$, then if $\ell +1$ is the coordinate for the first 
$h$ after $k$ we  have that  
$$
x = b_1 + \cdots + b_{\ell} \in r_{\Omega_n}(w)\setminus u,
$$
so $\Omega_n$ is a pre-shelling.   
\end{proof}     
If we define $\MAJ_{\ell} : \Dy_n \rightarrow \N$ by 
$$
\MAJ_{\ell}(w) = \sum_{i \in \Long(w)}i,
$$
we now have:
\begin{cor}
For all $n \geq 1$ we have 
$$
\beta_n(S) = |\{ w \in \Dy_n : \Long(w) = S \}|, 
$$
In particular the bi-statistic $(\lnfs,\MAJ_{\ell})$ has the 
$q$-Narayana distribution.
\end{cor} 
The Narayana statistic $\ea$ cannot in a natural way be associated to 
a shelling of $\Delta(J( {\bf 2} \times {\bf n}))$. However, 
it would be interesting to find a co-statistic $s$ for $\ea$ such that 
the bi-statistic $(\ea, s)$ has the $q$-Narayana distribution. 
%
%
\bibliographystyle{plain}

\end{document}